\documentclass[12pt]{amsart}
\include{amsmath} \include{amssymb} 
\usepackage{mathrsfs, fullpage}

\allowdisplaybreaks
\newcommand{\bB}{{\bf B}}
\newcommand{\bP}{{\bf P}}
\newcommand{\bR}{{\bf R}}
\newcommand{\bZ}{{\bf Z}}

\newcommand{\rE}{{\rm E}}
\newcommand{\scrI}{{\mathscr I}}
\newcommand{\coloneq}{\mathrel{\mathop:}=}

 \newtheorem{theorem}{Theorem}[section]

\title{Discovering An Algebra of Classes in the Algebra of Numbers---from 
George Boole to the Present}

\author{Stanley Burris}
\email{snburris@uwaterloo.ca}

\date{\today}

\begin{document}
\maketitle 

Many people mistakenly believe Boole's 1854 classic \cite{Boole-1854}, {\em The Laws of Thought}, 
must be about his creation of the subject called Boolean Algebra. 
Actually he used the symbols and equational reasoning of the ordinary number 
system, augmented by requiring that variables be idempotent, to develop his 
Algebra of Logic (AoL). 
Trying to understand his AoL as either Boolean Algebra or Boolean Rings is a 
mistake---nonetheless some expressions in Boole's AoL can easily be read as
belonging to these two subjects. 

A down-side of Boole's Algebra of Numbers approach was, after choosing
multiplication to be the totally defined operation of intersection, his operations 
of addition and subtraction on classes were {\em forced} to be partially defined---yet 
he carried out algebraic manipulations as though they were totally defined. 
His principle that
one could derive meaningful results about classes by methods employing 
uninterpretable terms was not well-received---others soon replaced 
Boole's AoL with what would eventually become Boolean Algebra.\footnote
{In 1864 Jevons \cite{Jevons-1864} scrapped Boole's numerical superstructure and 
replaced {\em classes} with {\em properties}. 
Converting his properties back into classes one finds that, starting with the 
natural operations of union, intersection and complement, he launched the development 
of modern Boolean Algebra, followed by such luminaries as Peirce, Venn, Schr\"oder, 
Huntington, Royce and Stone.
This note, however, is concerned with results that continued along the lines of Boole's 
Algebra of Numbers approach.} 
Although his AoL has more in common with Boolean Rings than 
Boolean Algebra,\footnote
{All the quasi-equations that hold in the Algebra of Numbers hold in Boole's AoL.
The equations also hold in Boolean Rings, but not all the quasi-equations---for 
example, Boolean Rings are not torsion-free.}
it was only much later that Boolean Rings appeared---taking the symmetric
difference as a fundamental operation was not an obvious choice.

Section \ref{sec Boole} is a summary of Boole's AoL. Section \ref{sec Whitney} is 
Whitney's 1933 encoding
of the modern Boolean Algebra of classes into the Algebra of Numbers using 
characteristic functions. 
Whitney did not realize how deeply connected his work was with that of Boole---such 
would have to wait till Hailperin's book of 1976/1986 with its focus on signed multisets 
(see Section \ref{sec Hailperin}). 
Section \ref{sec R01} explains Boole's Rule of 0 and 1 and the recent role of Horn 
sentences in proving his results. Section \ref{sec Brown} looks at Brown's
version of Boole's AoL using an algebra of proto-Boolean polynomials.

To achieve a smooth narrative the various original notations of the authors mentioned will 
frequently be set aside. 

\section{George Boole's Algebra of Logic} \label{sec Boole}

In \cite{Boole-1854}, p.~27, Boole asserted that only a few symbols were needed to 
develop an Algebra of Logic:
\begin{quote}
\textit{All the operations of Language, as an instrument of reasoning, 
may be conducted by a system of signs composed of the following elements,
viz.:}

1st. \textit{Literal symbols, as $x$, $y$, \&c., representing 
things {as subjects  of our conceptions}.} 

2nd. \textit{Signs of operation, as $+$, $-$, $\times$, standing for those 
operations {of the mind} by which {the conceptions of} things are combined 
or resolved
{so as to form new conceptions involving the same elements}.}

3rd. \textit{The sign of identity, $=$.}

\textit{And these symbols of Logic are in their use subject to definite 
laws, partly agreeing with and partly differing from the laws of the 
corresponding symbols in the science of Algebra.} 
\end{quote}
He would soon add the symbols 1 and 0 to the above.
Boole did not
have the following modern recursive definition of {\em terms} as strings of symbols 
for his AoL:
\begin{itemize}
\item Variables are terms, as are the symbols 0 and 1.
\item Suppose $s$ and $t$ are terms. Then so are $(s\cdot t)$, $(s+t)$, and $(s-t)$.
\end{itemize}
His notation followed that of the practice in the Algebra of Numbers, an informal 
mixture of terms and polynomials.

Given a universe of discourse, his model was to let 1 denote the universe, 0 the empty
class, and to define the three binary operations for classes $x$ 
and $y$ contained in the universe as follows:\footnote
{Modern notation is used on the right sides of the equal signs, where ``$\coloneq$'' 
means {\em equals by definition}.}
\begin{itemize}
\item
$x\times y \coloneq  x \cap y$
\item
$x+y \coloneq  \bigg\{
\begin{array} {l}
x \cup y \quad \text{if }\ x \cap y = 0\\
\text{undefined otherwise}
\end{array}$
\item
$x-y \coloneq  \bigg\{
\begin{array} {l}
x \setminus y \quad \text{if }\ y \subseteq x\\
\text{undefined otherwise.}
\end{array}$
\end{itemize} 
Instead of `undefined' Boole said `uninterpretable'.
In the following the product of $x$ and $y$ will be written either $x\cdot y$, or as 
Boole often 
wrote, simply $xy$ (just as one does in the algebra of numbers).
Clearly $1-x$ is the complement of the class $x$. The {\em idempotent} law $x^2=x$ 
for classes followed from his definition of multiplication. A term $t(x_1,\ldots,x_n)$
that is interpretable for all choices of classes for the $x_i$ is {\em totally} interpretable.\footnote
{One can easily give a recursive definition of the collection of totally interpretable terms,
and one can show that the domain of interpretability of a term $t(x_1,\ldots,x_n)$ is given
by the solutions of some equation $d(x_1,\ldots,x_n)=0$ where $d$ is totally interpretable.}
An equation $s=t$ is totally interpretable if both $s$ and $t$ are totally interpretable.

Boole expressed the standard operations on classes by totally interpretable terms  
(thanks to ``$-$'' being binary subtraction)
as well as by polynomials in his Algebra of Logic:\footnote
{Modern notation is used on the left sides of the equal signs.}
\begin{alignat*}{2}
{x\cap y}\  &=\  xy &&\\
{x\cup y}\  &=\ x + (1-x)y &\ &=\  x + y - xy\\
{x\triangle y}\  &=\ x(1-y) + (1-x)y\  &\  &=\ x + y - 2xy \\
{x'}\  &=\  1 - {x} &&\\
{U}\  &=\  1&&\\
{\O}\  &=\ 0.&&
\end{alignat*}
He used this to express propositions about classes as equations.
One can use Boole's AoL to derive the laws of Boolean algebra (but Boole
did not do this),  
for example, 
to prove the absorption law $x \cap (x \cup y) = x$ one has
\begin{eqnarray*}
x \cap (x \cup y) &=& x(x + y - xy)\\
&=& x^2 + xy - x^2y\\
&=& x + xy - xy \ = \ x.
\end{eqnarray*}

The Boolean Algebra and Boolean Ring of classes exist just below the surface of 
Boole's AoL---any totally interpretable term is equivalent to the Boolean Algebra 
term obtained by replacing addition
 by union, multiplication by intersection, and each subterm $s-t$ by $s\cap t'$.
But this is not how Boolean Algebra came into existence---see the earlier footnote 
about Jevons.
Any totally interpretable term can also be viewed as a Boolean Ring term. 
The Boolean Algebra and Boolean Ring of classes
are embedded in Boole's AoL, but the converse is false.

\subsection{The Rule of 0 and 1}

Boole introduced a powerful method to determine the valid equations and equational 
arguments in his algebra of logic on p.~37 of \cite{Boole-1854}, namely it suffices to check
them in the algebra of numbers when the variables are restricted to 0 and 1:
\begin{quote}
{\sf Let us conceive, then, of an Algebra in 
which the symbols $x$, $y$, $z$, etc. admit indifferently of the values 
$0$ and $1$, and of these values alone. The laws, the axioms, and 
the processes, of such an Algebra will be identical in their whole 
extent with the laws, the axioms, and the processes of an Algebra
of Logic. Difference of interpretation will alone divide 
them. Upon this principle the method of the following work is 
established. } 
\end{quote}
The strength of this principle was not understood for a century and a half---it 
will be discussed in Section \ref{sec R01}.
\subsection{Constituents}
Boole based the development of his algebra on the use of {\em constituents}.
The constituents of $x$ are simply $x$ and $1-x$; for $x,y$ they are
$xy$, $x(1-y)$, $(1-x)y$ and $(1-x)(1-y)$. 
In general for variables
$\vec{x} \coloneq  x_1,\ldots,x_m$ the constituents are the $C_\sigma(\vec{x})$ given
by
$$ C_\sigma(\vec{x})\ \coloneq\ \prod_{i=1}^m C_{\sigma_i}(x_i) $$
where $\sigma$ is a string of 0s and 1s of length $m$, and
$$ C_{\sigma_i}(x_i) \coloneq \bigg\{
\begin{array}{l l}
x_i & \text{if }\ \sigma_i = 1\\
1- x_i & \text{if }\ \sigma_i = 0.
\end{array}
$$
Boole noted that the $C_\sigma(\vec{x})$ are idempotent, pairwise disjoint and
their sum is 1.

\subsection{Development Theorem}
Perhaps the single most important result for his algebra of logic was that for any 
term $t(\vec{x},\vec{y})$, where $\vec{y} \coloneq  y_1,\ldots,y_n$,
$$ t(\vec{x},\vec{y}) = \sum_\sigma t(\sigma,\vec{y})\cdot C_\sigma(\vec{x}).$$
When there are no $y_i$ then one has the {\em complete} development
$$ t(\vec{x}) = \sum_\sigma t(\sigma)\cdot C_\sigma(\vec{x})$$
with the coefficients $t(\sigma)$ being integers.
The right side of this equation is totally interpretable iff each $t(\sigma)$ is idempotent,
that is, it is either 0 or 1. 
This is precisely the condition for $t(\vec{x})$ to be 
idempotent, that is, $ t(\vec{x})^2 = t(\vec{x})$.
The {\em idempotence} of a term is Boole's {\em condition of
interpretability}---a term is equivalent to a totally interpretable term iff it is idempotent. 
Note that $x+y$ does not satisfy this condition as $1+1$ is
neither $0$ nor $1$; $x+y$ is only interpretable under the condition that $xy = 0$.

 Boole stated that two terms $s(\vec{x}), t(\vec{x})$ were equal in his algebra iff 
 they had the same complete developments, that is, $s(\sigma) = t(\sigma)$ for
 all $\sigma$.
 In Boole's Algebra every equation $t_1(\vec{x}) = t_2(\vec{x})$ is equivalent to
 an equation $t(\vec{x}) = 0$, namely put $t(\vec{x}) \coloneq t_1(\vec{x}) - t_2(\vec{x})$.
 The Development Theorem, the properties of constituents, and the torsion-free property
  show that an equation $t(\vec{x}) = 0$ is equivalent
 to the collection of {\em constituent equations} $C_\sigma(\vec{x}) = 0$ where
 $t(\sigma) \neq 0$. 
 Thus any equation $t(\vec{x})=0$ is equivalent to an equation 
 $t^\star(\vec{x})=0$ that is totally interpretable, namely  let 
 $$t^\star(\vec{x}) \coloneq  \sum \{C_\sigma(\vec{x}) : t(\sigma)\neq 0\}.$$
  
\subsection{Reduction Theorem}

This says that any collection of equations 
$$t_1(\vec{x}) = 0,\ldots, t_k(\vec{x}) = 0$$
can be reduced to the single equation
$$t_1(\vec{x})^2 +\cdots+ t_k(\vec{x})^2\ =\ 0.$$
This follows from the torsion-free property and properties of constituents---they
imply that the collection of equations is equivalent to the collection
of constituent equations $C_\sigma(\vec{x}) = 0$ for which some $t_i(\sigma)\neq 0$.
Clearly $t_i(\sigma)\neq 0$ for some $i$  iff $t_1(\sigma)^2 +\cdots+ t_k(\sigma)^2\ \neq\ 0.$

\subsection{Elimination Theorem}

This says that the complete result of eliminating the
variables $\vec{x}$ from an equation $t(\vec{x},\vec{y})=0$ is the equation
$$ \prod_\sigma t(\sigma,\vec{y})\ = \ 0.$$
Thus eliminating $y$ from $t(\vec{x},y)=0$ gives $t(\vec{x},0)\cdot t(\vec{x},1) = 0$.
Boole did not define what he meant by the {\em complete result}.

\subsection{Solution Theorem}

Boole claimed that the solution for $y$ in an equation $t(\vec{x},y)=0$ consisted
of a pair of equations
\begin{eqnarray}
t(\vec{x},0)\cdot t(\vec{x},1) &=&0\\
y &=& t^\star(\vec{x},0) + v\cdot \big(1- t^\star(\vec{x},0)\big)\cdot \big(1- t^\star(\vec{x},1)\big),
\quad \text{for some } v.
\end{eqnarray}
(This is a tidied-up version of Boole's solution.) 
Boole frequently used this result in his examples in \cite{Boole-1854}, and for the application of
his AoL to Probability Theory.

\section{Hassler Whitney encodes Modern Boolean Algebra in the Algebra of Numbers 
using Characteristic Functions}\label{sec Whitney}

In 1933 Whitney \cite{Whitney-1933} showed that one did not need to be familiar with Boolean Algebra\footnote
{Like many, Whitney still used the traditional name {\em Algebra of Logic} instead of 
{\em Boolean Algebra}, the latter a name which had been promoted by the philosopher
 Josiah Royce of Harvard and his students.}
to verify
equations and equational arguments in this algebra---one could convert them into equations 
and equational arguments in the Algebra of Numbers. 
His main tool was {\em characteristic functions}
$\chi_A$ mapping a universe $U$ into $\{0,1\}$, taking the value 1 on $A$, the
value 0 off $A$.\footnote
{Evidently using characteristic functions was considered
novel as Whitney gave a reference to a 1916 work of de la Vall\'ee Poussin for the definition.}

The basic idea was that the characteristic function of a Boolean combination of
$A_i$ could be expressed as a polynomial in the characteristic functions of the
$A_i$, a polynomial with integer coefficients, since
\begin{eqnarray*}
\chi_{A\cap B} &=& \chi_{A}\cdot \chi_{B}\\
\chi_{A\cup B} &=& \chi_{A} + \chi_{B} - \chi_{A}\cdot \chi_{B}\\
\chi_{A\triangle B} &=& \chi_{A} + \chi_{B} - 2\chi_{A}\cdot \chi_{B}\\
\chi_{A'} &=& 1 - \chi_{A}\\
\chi_{U} &=& 1, \quad 
\chi_{\O} \ =\ 0.
\end{eqnarray*}
A polynomial in characteristic functions defines a function in the ring $\bZ^U$, 
but Whitney did not mention this ring.
He did not see the connection with Boole's AoL, namely if one restricted 
the ring operations to the characteristic functions then one had a partial algebra 
isomorphic to Boole's partial algebra.\footnote
{The restriction needs to use the binary subtraction operation, not the unary
minus operation.}
Although he had Boole's Development Theorem (which he called his First Normal Form)
 and the reduction of a system of equations
to constituent equations, there was no Elimination Theorem or Solution Theorem.

Whitney briefly considered arbitrary functions from the universe $U$ to the 
integers, generalizing the polynomials in the characteristic functions, but did 
not make explicit that these were precisely the elements of the ring $\bZ^U$.

\section{Theodore Hailperin's Axioms and Signed Multisets} \label{sec Hailperin}

Hailperin \cite{Hailperin-1976} noted that by identifying subsets of the universe $U$ 
with their characteristic functions, the ring $\bZ^U$ formed an extension
of Boole's partial algebra that made all terms in Boole's AoL interpretable.
 $\bZ^U$ is a non-trivial torsion-free\footnote
 {{\em Torsion-free} means that $nx = 0$ implies $x=0$, for $n=1,2,\ldots$}
 commutative ring with unity---Hailperin called members 
 of this ring {\em signed multi-sets}.  
He axiomatized Boole's AoL by the first-order axioms of non-trivial torsion-free
commutative rings with unity,\footnote
{Perhaps Hailperin was the first to recognize the critical role that the
torsion-free property played in the construction of Boole's AoL. 
Actually Boole only used it to show that $nC_\sigma(\vec{x}) = 0$ implies 
$C_\sigma(\vec{x}) = 0$. 
By the Development Theorem this leads to $nt(\vec{x}) = 0$ implies 
$t(\vec{x}) = 0$ for any term $t(\vec{x})$.
In an unpublished note Boole used the torsion-free property to prove that
$(x+y)^2 = x+y$ implies $xy = 0$, forcing him to define addition only for
disjoint classes.
Hailpern also included the unnecessary axiom $x^2 = 0$ implies $x=0$,
i.e., there are no nilpotent elements.} 
and proceeded to prove each of Boole's main theorems---expressed as a first-order
sentence  $\varphi$---by proving $\varphi_E$ from these axioms, where $\varphi_E$ 
is the result of restricting the quantifiers of $\varphi$ to idempotent elements. 
Finally Boole's work on logic had a firm foundation.
Note that Boole's rather awkward idempotent law, which applied only to variables, was 
dropped from the axioms by Hailperin---it was implicitly replaced by the restriction of the 
quantifiers to idempotent elements.

\section{R01 and Horn sentences} \label{sec R01}

Boole's Rule of 0 and 1, abbreviated as R01, was finally given its proper formulation in 2013 
in \cite{BuSa-2013}. Following Hailperin's lead, a first-order sentence $\varphi$ will be said 
to hold in Boole's AoL if
$\varphi_E$ follows from the above-mentioned axioms of Hailperin.\footnote
{Boole surely viewed his models as being {\em atomic}, that is, each non-empty class contains
an atom (a singleton). So to extend Boole's theorems to first-order properties $\varphi$ not
considered by him, one needs to add an axiom expressing the atomic property to Hailperin's axioms.}
By the completeness
theorem for first-order logic this is equivalent to saying that $\varphi_E$ holds in all rings
satisfying Hailperin's axioms. Let $\bR$ be such a ring, and let $\bR_E$ be the subring generated
by the idempotents of $\bR$.\footnote
{For the ring $\bZ^U$ the subring generated by the idempotents consists of all members of 
$Z^U$ with a finite range.}
Clearly $\bR \models \varphi_E$ iff $\bR_E \models \varphi_E$.
It is known that $\bR_E$ is isomorphic to a bounded Boolean power $\bZ[\bB]^\star$ of
$\bZ$ (see  \cite{Keimel-1970}, \cite{Subramanian-1969}). 
From \cite{BuWe-1979} we know that Horn sentences true in $\bZ$
also hold in $\bZ[\bB]^\star$, giving the following result.
\begin{theorem} \label{R01}
A Horn sentence $\varphi$ holds in Boole's Algebra of Logic
iff\; $\bZ \models \varphi_\rE$.
\end{theorem}
Note that $\bZ \models \varphi_\rE$ simply says that $\varphi$ holds in
the ring of integers when the variables are restricted to the values 0 and 1.\footnote
{Thus Boole's AoL satisfies all the equations and quasi-equations
that hold for the integers.
Some of the more exciting quasi-equations that hold in the integers become quite
trivial in Boole's AoL, for example, $x^n + y^n = z^n \rightarrow xyz = 0$ for $n\ge 3$
becomes $x+y=z \rightarrow xyz = 0$ in Boole's AoL, a quasi-equation that 
 R01 immediately validates.}
Since each of Boole's main theorems can be expressed as a Horn sentence $\varphi$,
Theorem \ref{R01} readily proves that they hold in Boole's AoL.

\section{Frank Brown's Proto-Boolean Algebras} \label{sec Brown}

In his 2009 paper \cite{Brown-2009} Brown said that Boole did not base his AoL on
signed multisets and rings but on
the algebra of polynomials, and he set out to give an updated presentation of the latter. 
Brown worked with polynomials in which variables can only occur to the first power. 
He called these {\em proto-Boolean polynomials}, abbreviating the name to 
$p$-{\em polynomials}. The sum and difference of two $p$-polynomials is again
a $p$-polynomial, but not in general their product. A new multiplication $q \odot r$
is introduced for $p$-polynomials that takes the ordinary product $q\cdot r$ and
reduces to 1 all exponents greater than 1. For example, 
$xy\odot (x+y) = 2xy$.
The new product of $p$-polynomials is again a $p$-polynomial.

Let $P(X)$ be the set of $p$-polynomials with variables from $X$.\footnote
{Brown only considered finite $X$.}
With the usual plus and minus operations along with the new multiplication one has 
a ring $\bP(X)$ that satisfies Hailperin's axioms---indeed it is isomorphic 
to the quotient ring $\bZ[X]/\scrI$ where 
$\scrI$ is the ideal of $\bZ[X]$ generated by $\{x^2-x : x \in X\}$. 
The isomorphism
is given by mapping each $p$-polynomial $q$ to the coset $q+\scrI$.
This ring in turn is isomorphic to $\bZ^{2^n}$ if $|X|=n$.

Brown proved Boole's results hold in these rings of $p$-polynomials by proving
the appropriate $\varphi_E$ where the variables range over the idempotent
polynomials in $P(X)$. 
The connection with theorems about classes is not as obvious as when working
with Hailperin's models $\bZ^{U}$ and letting the variables range over the 
 idempotent elements, namely over the characteristic functions. 
Brown's connection was simply to point out that, for $|X|=n$, the Boolean Algebra 
of idempotents in the ring $\bP(X)$ is isomorphic to a Boolean algebra of $2^{2^n}$ 
sets freely generated by $n$ sets.
\smallskip

{\sc References}   \label{references}


\begin{enumerate}

\bibitem{Boole-1854}
George Boole,
{\em An Investigation of The Laws of Thought on Which are Founded the
Mathematical Theories of Logic and Probabilities}.
 Originally published by Macmillan,
London, 1854. Reprint by Dover, 1958.

\bibitem{Brown-2009}
Frank Markham Brown,
{\em George Boole's deductive system}.
Notre Dame J.~Form.~Log., {\bf 50} (2009), 303--330.

\bibitem{Burris-2015c}
Stanley Burris,
{\em George Boole and Boolean Algebra}. Eur.~Math.~Soc. Newsletter, Dec 2015, 27--31.
 
 \bibitem{BuWe-1979}
 \bysame and Heinrich Werner,
 {\em Sheaf constructions and their elementary properties.}
 Trans.~Amer.~Math.~Soc., {\bf 248}, No. 2 (1979), 269-309.
 
    \bibitem{BuSa-2013}
\bysame and H.P.~Sankappanavar,
 {\em The Horn Theory of Boole's Partial Algebras}. Bull. Symb.~Log., 
 {\bf 19} (2013), 97--105.
 
 \bibitem{Hailperin-1976}
 Theodore Hailperin,
{\em Boole's Logic and Probability}, Series: Studies in Logic and the Foundations of 
Mathematics, {\bf 85}, Amsterdam, New York, Oxford: Elsevier North-Holland. 1976. 
2nd edition, Revised and enlarged, 1986.

 \bibitem{Hailperin-1981}
\bysame,
{\em Boole's Algebra isn't Boolean algebra}. Math.~Mag., {\bf 54}, 
No. 4 (1981),  172--184.
 
 \bibitem{Jevons-1864}
William Stanley Jevons, {\em Pure Logic, or the Logic of Quality apart from Quantity: with
Remarks on Boole's System and on the Relation of Logic and Mathematics}.  Edward
Stanford, London, 1864. Reprinted 1971 in Pure Logic and Other Minor Works ed. by R.
Adamson and H.A. Jevons, Lennox Hill Pub. \& Dist. Co., NY.

\bibitem{Keimel-1970}
Klaus Keimel, {\em Alg\`ebres commutatives engendrees par leurs elements idempotents.}
Can.~J.~Math., {\bf XXII}, No. 5 (1970), 1071--1078.

\bibitem{Subramanian-1969}
H. Subramanian, 
{\em Integer-valued continuous functions.} Bull.~Soc.~Math.~France, {\bf 97} (1969),
275--283.

\bibitem{Whitney-1933}
Hassler Whitney,
{\em Characteristic functions and the algebra of logic.}
Ann.~of Math., {\bf 34} (1933), 405--414.

\end{enumerate}
\end{document}